\documentstyle[12pt]{article}
\newcommand{\be}{\begin{equation}}
      \newcommand{\ee}{\end{equation}}
      \newcommand{\ba}{\begin{eqnarray}}
       \newcommand{\ea}{\end{eqnarray}}
\newcommand{\ban}{\begin{eqnarray*}}
       \newcommand{\ean}{\end{eqnarray*}}

\newcommand{\lp}{\langle}
\newcommand{\rp}{\rangle}
\newcommand{\ra}{\rightarrow}

 \newcommand{\qed}{\hspace*{\fill}\rule{3mm}{3mm}\quad}
 \newcommand{\Pf}{\noindent {\em Proof.} }

\newcommand{\Rk}{\noindent {\em Remark} }

\newcommand{\sect}[1]{\section{#1} \setcounter{equation}{0}}

\newtheorem{theo}{Theorem}[section]

\begin{document}
\newtheorem{lem}[theo]{Lemma}
\newtheorem{prop}[theo]{Proposition}
\newtheorem{coro}[theo]{Corollary}

\title{Controlled Geometry via Smoothing  \footnote{ 1991 {\em Mathematics
Subject Classification}. Primary 53C20.}}
\author{Peter Petersen\thanks{Partially supported by NSF and NYI grants} \\
  \and Guofang Wei\thanks {Partially supported by NSF Grant \# DMS9409166.
}\\
\and Rugang Ye\thanks{Partially supported by NSF Grant \# DMS9401106}
}
\date{}
\maketitle

\begin{abstract}
We prove that Riemannian metrics with a uniform weak norm can be smoothed to
having arbitrarily high regularity. This generalizes all previous smoothing
results. As a consequence we obtain a generalization of Gromov's almost flat
manifold theorem. A uniform Betti number estimate is also obtained.
\end{abstract}

\newcommand{\Tr}{\mbox{Tr}}
\newcommand{\inj}{\mbox{inj}}
\newcommand{\vol}{\mbox{vol}}
\newcommand{\diam}{\mbox{diam}}
\newcommand{\Ric}{\mbox{Ric}}
\newcommand{\conj}{\mbox{conj}}
\newcommand{\Hess}{\mbox{Hess}}
\newcommand{\divg}{\mbox{div}}
\sect{Introduction}

An ultimate goal in geometry is to achieve a classification scheme, using
natural geometric quantities to characterize the topological type
or diffeomorphism type of Riemannian manifolds. While this grand scheme  seems
to be an impossible dream, its basic philosophy has been
 a driving force in many important developments in Riemannian geometry. The
sphere theorems and various topological finiteness theorems
 are typical examples. These results are concerned with control of global
topology of manifolds, and a crucial point therein is to control,
 uniformly, the local topology.

Control of local topology often follows from control of local
 geometry. Here, by local geometry, we mean the local behavior of the metric
tensor. On the other hand, control of local geometry is frequently also
 the essential ingredient for control of global geometry, such as in
 Cheeger-Gromov's compactness theorem and its various extensions,
which
can be named geometric finiteness theorems. Notice that some rudimental
topological finiteness results are direct corollaries of
geometric
finiteness theorems. But the significance of the latter  goes beyond this. In
any case, control of local geometry is obviously a key
topic. An interesting and important aspect of this topic is various degrees of
control of local geometry needed or available in different
situations. In \cite{p}, the first author introduced a sequence of norms which
provide a certain quantitative measure for local geometric
control.  These norms can be defined either in terms of the
$C^{k,\alpha}$-norms
 or the $L^{k,p}$-norms for functions, and they are defined on a given
 scale. For example, the $C^{k,\alpha}$-norm of a Riemannian manifold on scale
$r$ is bounded, if it is covered by coordinate charts of size comparable to $r$
such that the metric tensor expressed in the coordinates
is uniformly bounded in $C^{k, \alpha}$-norm, and that the coordinate
transition functions are uniformly bounded
in $C^{k+1, \alpha}$-norm. Note that the local topology is uniformly trivial if
one of these norms on some scale is bounded. To admit richer topological and
geometric structures under norm bounds, we shall introduce a weak version
of these norms. The essential new  feature is that we allow coordinate maps
to have double points. In spirit, this is similar to replacing an injectivity
radius bound by a conjugate radius bound. (Of course, e.g. a weak (harmonic)
$C^{0, \alpha}$ bound is so weak that it is far from implying a conjugate
radius bound.)

Indeed, our basic theme is to try to find the minimal degree of control of
local geometry under which interesting geometric and
topological
consequences can be
drawn. Traditional geometric conditions such as curvature bounds imply various
degrees of local geometric control, so we can think from their perspectives.
Historically,  sectional curvature bounds were the first to be systematically
studied. They can roughly be compared with weak  $C^2$-norm bounds, at least
the latter imply
the former. Since understanding of sectional curvature bounds  has  been
reached on a
good level, it is natural to try to find the minimal degree of local geometric
control under which a metric can be approximated by metrics with  sectional
curvature bounds or weak $C^2$-norm bounds (or better bounds). In this paper,
we present a result towards this goal,
along with some applications.  The  local geometric control we need is as weak
as  a bound on the weak harmonic $C^{0, \alpha}$-norm,  or a bound on the weak
$L^{1, p}$-norm. These do appear to be the sought-after minimal
degree of local geometric control in our set-up.

To formulate the result precisely, we introduce the following classes of
Riemannian manifolds. (The definition of the weak norms are given in \S 2.)
\newline

\noindent {\bf Definition}. Given $n \geq 2$, $0 < \alpha < 1$, $p > n$ and
function $Q: (0, \infty) \ra [0, \infty)$ which is nondecreasing in $r$ and
satisfies $\lim_{r \ra 0} Q(r) = 0$, we define
\[
{\cal M}(n, \alpha, Q) = \left\{ (M, g) \left| \begin{array}{l} (M,g) \mbox{ is
a complete Riemannian manifold}, \dim M = n, \\
\mbox{the weak harmonic}\ C^{0, \alpha}\ \mbox{norm}\
\|(M,g)\|^{W,h}_{C^{0,\alpha}, r} \leq Q(r)\\ \mbox{ for all positive }\ r \leq
1
\end{array} \right.
\right\},
\]
and
\[
{\cal M}(n, p, Q) = \left\{ (M, g) \left| \begin{array}{l} (M,g) \mbox{ is a
complete Riemannian manifold}, \dim M = n, \\
\mbox{the weak}\ L^{1,p}\ \mbox{norm}\ \|(M,g)\|^{W}_{L^{1,p}, r} \leq Q(r)\\
\mbox{ for all positive}\ r\leq 1
\end{array} \right.
\right\}.
\]

\Rk 1. By Anderson-Cheeger's work \cite{ac} manifolds with a lower bound for
Ricci curvature and a positive lower bound for conjugate radius belong to these
two classes.

\Rk 2. It is unknown whether (weak) $L^{1,p}$ bounds imply (weak) harmonic
$C^{0, \alpha}$ bounds. (By the Sobole
v embedding, they do imply $C^{0,\alpha}$ bounds, but not yet harmonic $C^{0,
\alpha}$ bounds. ) In other words,
it is unknown whether controlled $C^{0,\alpha}$  harmonic coordinates exist
under the assumption of a (weak) $L^{
1, p}$ bound.  This question seems rather subtle. We plan to return to it in
the future. At the moment, we consider the said two bounds as independent
conditions. Note that one can further ask  whether (weak) $C^{0, \alpha}$
bounds
imply (weak) harmonic $C^{0, \alpha}$ bounds. There does not seem to be any
evidence to support an answer in the affirmative.
 \\

\begin{theo} \label{main}
For every manifold $(M,g)$ in ${\cal M}(n, \alpha, Q)$ and positive numbers
$\epsilon $,   there  are  metrics $g_\epsilon$ on $M$ such that
\ba
e^{-\epsilon}g   \leq &  g_\epsilon &\leq e^\epsilon g  \label{t1},\\
\| (M,g_\epsilon)\|^{W}_{C^{0,\alpha}, r} & \leq & 2Q(r) \label{t2}, \\
\| (M,g_\epsilon)\|^W_{C^{k,\alpha}, r} & \leq & \widetilde{Q},
\ea
\noindent where $k$ is an arbitary positive integer, $0 < r \leq 1$ and
$\widetilde{Q} = \widetilde{Q} (n, k , \epsilon, \alpha, Q(r))$
denotes a positive number depending only on $n, k, \epsilon, \alpha$ and
$Q(r)$.
\end{theo}

\begin{theo} \label{main'}
For every  closed manifold $(M,g)$ in ${\cal M}(n, p, Q)$ and positive numbers
$\epsilon $,   there are  metrics $g_\epsilon$ on $M$ such that
\ba
e^{-\epsilon}g   \leq &  g_\epsilon &\leq e^\epsilon g  \label{t1'},\\
\| (M,g_\epsilon)\|^{W}_{L^{1,p}, r} & \leq & 2Q(r) \label{t2'}, \\
\| (M,g_\epsilon)\|^W_{L^{k,p}, r} & \leq & \widetilde{Q},
\ea
\noindent where $k$ is an arbitary positive integer,  $0 < r \leq 1$ and
$\widetilde{Q} = \widetilde{Q} (n, k , \epsilon, p, Q(r))$  denotes
 a positive number depending only on $n, k, \epsilon, p$ and $Q(r)$.
\end{theo}

\Rk  Theorem 1.2 actually also holds for complete, non-compact manifolds. This
will be  shown in a paper by the third author. \\

Thus a metric with some regularity (given by the weak norm) can be deformed or
smoothed to a nearby one with arbitrarily high regularity. In particular, {\it
manifolds with a lower bound on Ricci curvature and a positive lower bound on
conjugate radius  can be smoothed.}  Previous smoothing results have been
concerned with metrics with various curvature bounds, and involved two
independent techniques: the embedding method and the Ricci flow. The embedding
technique in smoothing as
used by Cheeger-Gromov \cite{cg} consists of embedding (or immersing) a given
manifold into a Euclidean space and then perturbing it
suitably by a smoothing operation, which is based on the classical convolution
process. The smoothing result in \cite{cg} is that metrics on closed manifolds
with lower and upper bounds on sectional curvatures and  a positive lower bound
on injectivity radius can be smoothed to metrics with bounds on all derivatives
of the  Riemann curvature tensor. Later, by embedding into a Hilbert space
instead of a finite dimensional space, Abresch \cite{a} was able to remove  the
condition on injectivity radius and extend to complete manifolds. More recently
Shen \cite{sh} showed that manifolds with a lower bound on sectionl curvatures
and a positive lower bound on injectivity radius can be smoothed to having
two-sided sectional curvature bounds. The technique of Ricci flow is based on
the fundamental work of Hamilton \cite{ha}. Using this technique, Bemelmans-Min
Oo-Ruh \cite{bmr} obtained the same result as in \cite{cg} without injectivity
radius lower bound, and Shi \cite{s} obtained the same result as in \cite{a}.
Later work considers metrics with other kinds !
 of curvature bound.  For example,
in \cite{y1,y2}
Yang dealt with integral bounds on sectional curvatures. In \cite{dwy}, Ricci
curvature bounds were treated.

By virtue of the available constructions of controlled harmonic coordinates
 under various curvature bounds, all these smoothing results are consequences
of Theorem~\ref{main} or Theorem~\ref{main'}.

As typical applications we present the following two results.

\begin{theo}[Betti number estimate] For the class of manifolds $M^n$ in \\
${\cal M}(n, \alpha, Q)$ or in ${\cal M}(n, p, Q)$, and satisfying $\diam_M
\leq D$,
we have the estimate for the Betti numbers
\be
\sum_i b^i(M^n) \leq  C(n,D,\alpha,Q)\ \mbox{or}\ C(n,D,p,Q),\label{be}
\ee
and the estimate for the number of isomorphism classes of rational homotopy
groups
\be
 \pi_q(M) \otimes Q  \leq   C(n,q,D,\alpha,Q)\ \mbox{or}\ C(n,q,D,p,Q)\
\mbox{for}\ q \geq 2. \label{hge}
\ee
\end{theo}

(\ref{be}) follows from Theorem~\ref{main}, \ref{main'} and Gromov's uniform
betti number estimate regarding sectional curvature \cite{g2}. This estimate
can also be proved directly using Toponogov type comparison estimate introduced
in \cite{w}, see \cite{pw} for details. In \cite{w} the same estimate
(\ref{be}) is given for the class of manifolds satisfying $\Ric_M \geq
-(n-1)H$, conj $\geq r_0$ and $\diam_M \leq D$. (\ref{hge}) follows from
Theorem~\ref{main}, \ref{main'} and the results in  \cite{r}.

\begin{theo} \label{aflat}
There exists an $\epsilon =\epsilon(n,\alpha,Q)$ or $\epsilon(n,p,Q)> 0$ such
that if a manifold   $M^n$ belongs to
${\cal M}(n, \alpha, Q)$ or ${\cal M}(n, p, Q)$ and diam $\leq \epsilon$, then
$M$ is diffeomorphic to an infranilmanifold.
\end{theo}

This generalizes Gromov's almost flat manifold theorem \cite{g1} as well as its
 generalization in \cite{dwy}. (The proof is simple: combine Theorem 1.1 with
\cite{g1}.)

The proof of Theorem~\ref{main} and Theorem 1.2 uses the embedding method in
\cite{a}. Roughly speaking, we embed a given Riemannian manifold into
 the Hilbert space of $L^2$-functions on it, and then use the embedding to pull
back the $L^2$-metric of the Hilbert space.
The crucial point is of course to find a suitable embedding, such that the
pull-back metric will enjoy nice properties. In \cite{a}, the embedding is
defined in terms of distance functions. In our situation, these functions are
not appropriate, and we employ instead solutions of a canonical geometric
partial differential equation. Now if e.g. the  harmonic $C^{0,\alpha}$-norm of
the manifold is bounded, then a uniform pointwise bound on sectional curvatures
will hold for the
 pull-back metric, and hence we can apply the smoothing results for metrics
with sectional curvature bounds as given e.g. in \cite{a} or \cite{s}.

If we only assume that the weak harmonic $C^{0, \alpha}$-norm of the manifold
is bounded, i.e.  it is in the class ${\cal M}(n, \alpha, Q)$,
  the global embedding is generally not under control. To remedy the situation,
 we follow the idea in \cite{a} of employing instead local embeddings. In
\cite{a}, Abresch uses the exponential map to lift local patches of the
manifold and
his local embeddings are exactly  embeddings of these lifted patches.  In our
situation, the exponential map is not suitable. Our substitute for it is the
coordinate maps.
Thus we use them to lift local patches,  and construct embeddings of the lifted
patches via the same geometric partial differential equation as mentioned
before.
To make sure that the pull-back metrics induced by these local embeddings
descend to the local patches
and that the resulting metrics patch together to define a metric globally, it
is crucial to require the embeddings to be equivariant under isometries.  Since
our embeddings are defined in terms of solutions of  a canonical geometric PDE,
they
naturally share this equivariance property.

Basically, the above scheme also works for manifolds in the class ${\cal M}(n,
p, Q)$, but some modifications are necessary.
As before, the said pull-back metrics descend to
yield a new metric on the underlying manifold. But these metrics satisfy here
an integral bound on sectional curvatures rather than a poinwise bound. This is
a new situation. To handle it, we apply the Ricci flow and follow the arguments
in \cite{dwy}. A pointwise bound on Ricci curvature is used in several places
in\cite{dwy}. Since no such bound is available in our current situation, the
arguments in \cite{dwy} need to be improved and modified.  The result we thus
arrive at not only completes the smoothing scheme for the class ${\cal M}(n, p,
Q)$, but also provides some new understanding of short time existence of the
Ricci flow.



\sect{Norm, Weak Norm and Smoothing}
Fix an integer  $k \geq 0$ a number $0 \leq \alpha \leq 1$. The {\it
$C^{k,\alpha}$-norm } of an $n$-dimensional Riemannian manifold $(M,g)$ on
scale $r$, $\|(M,g)\|_{C^{k,\alpha}, r}$, is defined to be the infimum of
positive numbers $Q$ such that there exist embeddings: $$\varphi_s: B(0,r)
\subset R^n \ra U_s \subset M$$
($B(0,r)$ denotes the closed ball of radius $r$ centered at the origin) with
images $U_s$, $s \in  \cal S$ (an index set),  with the following properties:\\
1) $e^{-Q} \delta_{ij} \leq g_{s,ij} \leq e^Q \delta_{ij}$, \\
2) Every metric ball $B(p, \frac{r}{10}e^{-Q}),\ p\in M$ lies in some set
$U_s$,\\
3) $r^{|j|+ \alpha} \|\partial^j g_{s,ij}\|_{C^\alpha} \leq Q$ for all
multi-indices  $j$ with $0 \leq |j| \leq k$. \\
Here $g_{s,ij}$ denote the coefficients of $g_s=\varphi_s^*g$ on $B(0,r)$, and
$\delta_{ij}$ are the Kronecker symbols.

Note that this definition is slightly different from the corresponding one in
\cite{p}, where in addition the (rescaled) $C^{k+1, \alpha}$-norm of the
transition functions are required to be under control.  For convenience, we can
call the $C^{k, \alpha}$-norm (of Riemannian manifolds) as defined in \cite{p}
the {\it strong $C^{k, \alpha}$-norm}. (Note however that  the "strong"
harmonic $C^{k, \alpha}$-norm is equivalent to the harmonic $C^{k,
\alpha}$-norm.)

We define the harmonic $C^{k,\alpha}$-norm on scale $r$, $\|(M,g)\|^h_{C^{k,
\alpha}, r}$,   by requiring additionally the following\\
4) $\varphi_s^{-1}: U_s \ra R^n$ is harmonic, \\
which is equivalent to saying that \\
4$'$) id $: B(0,r) \ra B(0,r)$ is harmonic with respect to $g_{s}$ on the
domain and the Euclidean metric on the target, which is in turn equivalent to
saying
 that $$\sum_i \partial_i (g_s^{ij} \sqrt{ \det g_{s,ij}})= 0$$ for all $j$.

If $k \geq 1$ and $p >n$ (when $k=1$) or $p > \frac{n}{2}$ (when $k \geq 2$),
then we define the $L^{k,p}$-norm on the scale of $r$, $\|(M,g)\|_{L^{k,p},
r}$,  by retaining 1) and 2), and replacing 3) by\\
3) $r^{|j|-\frac{n}{p}} \|\partial^j g_{s,ij}\|_{L^p} \leq Q$ for all $1 \leq
|j| \leq k$.

The harmonic $L^{k,p}$-norm is defined similarly. For any choice of these
norms, it is clear that  the local topology is trivial on some uniform scale
for any
class of manifolds with uniformly bounded norm. (Note that the injectivity
radius may not be uniformly positive though.) To allow nontrivial local
topology, we introduce the  weak norms $\|\ \ \|^W_{C^{k,\alpha},r}$ and $\|\ \
\|^W_{L^{k,p},r}$, which are defined in identical ways except that each
$\varphi_s: B(0,r) \ra U_s$ is assumed to be  a {\it local} diffeomorphism
instead of diffeomorphism. The corresponding weak  harmonic norms $\|\ \
\|^{W,h}_{C^{k,\alpha},r}$ and $\| \ \ \|^{W,h}_{L^{k,p},r}$  are defined
in a similar way, with 4) being replaced by 4$'$).

Note that (weak) harmonic norms dominate (weak) norms on the same scale. We
also have $\|\ \ \|^W_{\ ,r} \leq \|\ \ \|_{\ ,r}$ and $\|\ \ \|^{W,h}_{\ ,r}
\leq \|\ \ \|^h_{\, r}$.    All norms are continuous and non-decreasing in $r$.
If $(M,g)$ is sufficiently smooth, these norms converge to zero as $r \ra 0$.
Furthermore, (weak) $C^{k,\alpha}\ (L^{k,p})$ norms vary continuously in the
$C^{k,\alpha}\ (L^{k,p})$ topology of Riemannian manifolds. See \cite{p} for
the relevant details.

We point out that $R^n$ is the only space with norm $=0$ on all scales. And
flat manifolds are the only spaces with weak norm $=0$ on all scales.

Conventional geometric conditions such as curvature bounds imply norm bounds.
Such implications are mostly contained in constructions of controlled harmonic
coordinates and are a crucial ingredient for  various compactness theorems. To
have a clear perspective, we collect these results in the following
proposition.
\begin{prop}
There is a $Q(H,i_0,r,p)$ with $\lim_{r \ra 0} Q(H,i_0,r,p) =0$ such that for
manifolds with \\
a) $|K| \leq H, \ \inj \geq i_0$, then $\|(M,g)\|^h_{L^{2,p},r} \leq
Q(H,i_0,r,p)$; \\
b) $|K| \leq H$, then  $\|(M,g)\|^{W,h}_{L^{2,p},r} \leq Q(H,r,p)$; \\
c) $|\Ric | \leq (n-1)H,\  \inj \geq i_0$, then $\|(M,g)\|^h_{L^{2,p},r} \leq
Q(H,i_0,r,p)$; \\
d) $\Ric \geq -(n-1)H,\ \inj \geq i_0$, then $\|(M,g)\|^h_{L^{1,p},r} \leq
Q(H,i_0,r,p)$; \\
e) $\Ric \geq -(n-1)H$, conj $ \geq i_0$, then $\|(M,g)\|^{W,h}_{L^{1,p},r}
\leq Q(H,i_0,r,p)$.
\end{prop}

These results follow from works of Jost-Karcher \cite{jk}, Anderson \cite{an}
and Anderson-Cheeger \cite{ac}.

We now turn to the smoothing question. As explained in the introduction, our
strategy is to first achieve sectional curvature bounds by embedding into the
Hilbert space of $L^2$-functions. This is done in the next two sections. The
higher regularity smoothing then easily follows from known smoothing results.
Consider $(M,g) \in {\cal M}(n, \alpha, Q)$. We have a collection of local
diffeomorphisms
\[
\varphi_s:\ B(0,r) \ra U_s \subset M
\]
satisfying 1), 2), 3) and 4$'$).

In the next section we will construct
a canonical embedding
\[
F_s: (B(0,r),g_s) \ra L^2 (B(0,r), g_s),
\]
 where $g_s = \varphi_s^*g$. We use $F_s$ to pullback the $L^2$ metric of $L^2
(B(0,r), g_s)$ to produce a new metric $\tilde{g}_s$ on $B(0,r)$. This
construction works for general metrics on $B(0,r)$, and has the following
equivariance property, which will be proved in \S 4. Namely, if $g_1,\ g_2$ are
two metrics on $B(0,r)$ such that there is an isometric embedding
\[
\psi: (B(0,r), g_1) \ra (B(0,r),g_2)
\]
and if $\tilde{g}_1,\ \tilde{g}_2$ are obtained via the above  construction,
then
\[
\psi: (B(0,r), \tilde{g}_1) \ra (B(0,r), \tilde{g}_2)
\]
is also an isometric embedding. Granted this (see Proposition~\ref{4.4}) we
have
\begin{prop}
There exists a smooth metric $\bar{g}$ on $M$ such that the pullback of
$\bar{g}$ by $\varphi_s$ is exactly  $\tilde{g}_s$.
\end{prop}
\Pf Let $r_1= \frac{r}{10} e^{-Q}$. Then for every $p \in M$, $B(p,r_1) \subset
U_s$ for some $s$. It follows that there exists a $\tilde{p} \in B(0,r)$ such
that $B_{g_s}(\tilde{p},r_1) \subset B(0,r)$ and $\varphi_s (\tilde{p}) = p$.

We now define the metric $\bar{g}$ as follows. If $X,Y \in T_pM$, then
\[
\bar{g}(X,Y) = \tilde{g}_s \left(((\varphi_s)_*|_{\tilde{p}})^{-1}(X),
((\varphi_s)_*|_{\tilde{p}})^{-1}(Y)\right).
\]

To show that this metric is well-defined, let $\tilde{p}'$ be another such
point, i.e. for some $s'$, $B_{g_{s'}}(\tilde{p}',r_1) \subset B(0,r)$ and
$\varphi_{s'} (\tilde{p}') =p$. Let $r_4 = \frac{r}{20} e^{-4Q}$ and $r_3 =
\frac{r}{20} e^{-3Q}$. Denote $g_0$ the Euclidean metric. Then we can show that
\begin{lem} \label{hom}
There is an isometric embedding
\[
\psi: (B_{g_0}(\tilde{p},r_4),g_s) \ra (B_{g_{s'}} (\tilde{p}',r_3),g_{s'}).
\]
\end{lem}
\Pf First $\psi$ can be defined as follows.
Since $g_s$ is $e^Q$-quasi-isometric to $g_0$,
\be  \label{r34}
B_{g_0}(\tilde{p}, r_4) \subset B_{g_s} (\tilde{p}, r_3).
\ee
 For any point $q \in B(\tilde{p},r_4)$, connect $q$ to the center point
$\tilde{p}$ with a curve $\tilde{\gamma}$ in $B_{g_0}(\tilde{p},r_4)$ such that
the length of $\tilde{\gamma}$ $l_{g_s}(\tilde{\gamma}) < r_3$. Since
\[
\varphi_s:\  (B_{g_s} (\tilde{p}, r_3), g_s) \ra B(p,r_3)
\]
is a local isometry and $\varphi_s (\tilde{p}) = p$. From (\ref{r34})
$\varphi_s$ maps the curve $\tilde{\gamma}$ to a curve $\gamma$ in $B(p,r_3)$
starting with $p$ and $l(\gamma) < r_3$. Again since $\varphi_{s'}$ is a local
isometry and $\varphi_{s'} (\tilde{p}') =p$. The curve $\gamma$ then can be
lifted via $\varphi_{s'}$ to a curve in $B_{g_{s'}}(\tilde{p}',r_3)$ starting
with $\tilde{p}'$. The other end point of this curve is defined to be the image
of $q$. (Note that, in general, lifting can not be done for incompelete space.
Here the map is a local isometry and the curve starts from the center, and we
have control on the length of the curve and the size of the metric ball, so it
will not hit the boundary during lifting.) Now we will show that $\psi$ is
well-defined, i.e. the image is independent of the choices of the curve
$\tilde{\gamma}$. If $\tilde{\gamma}_1$ is another curve in
$B_{g_0}(\tilde{p},r_4)$ connecting $q$ to the center point $\tilde{p}$ with
$l_{g_s!
 }(\tilde{\gamma}_1) < r_3$, we ca
$g_s$ is $e^Q$-quasi-isometric to $g_0$. Then $\varphi_s$ maps
$\tilde{H}(s,t)$ to a homotopy $H(s,t)$ in $B(p,2r_3)$ with $l(H(s,\cdot)) <
2r_3$ for each $s$. Therefore $H(s,t)$ can be lifted via $\varphi_{s'}$ to a
homotopy in $B_{g_{s'}}(\tilde{p}',2r_3)$  starting with $\tilde{p}'$. By the
(localized) homotopy lifting lemma the other end points are all the same.
Therefore $\psi$ is well-defined.

Next we show that $\psi$ is one-to-one. Let $r_2 =  \frac{r}{20} e^{-2Q}$.
Then
 $$B_{g_{s'}} (\tilde{p}',r_3) \subset B_{g_0}(\tilde{p}',r_2) \subset
B_{g_{s'}} (\tilde{p}',\frac{1}{2}r_1).$$ Since $B_{g_0}(\tilde{p}',r_2)$ is
 an Euclidean ball
one can construct ``inverse" $\phi$ similarly as above:
\[
\phi: \ (B_{g_{s'}} (\tilde{p}',r_3),g_{s'}) \ra
(B_{g_s}(\tilde{p},\frac{1}{2}r_1),g_s).
\]
Thus $\psi$ is one-to-one. That $\psi$ is an isometric embedding follows from
the construction.
\qed

Now using the equivariance, we have
\[
\psi^* \tilde{g}_{s'} = \tilde{g}_s.
\]
Therefore
\ban
\lefteqn{\tilde{g}_{s'}\left(((\varphi_{s'})_*|_{\tilde{p}'})^{-1}(X),
((\varphi_{s'})_*|_{\tilde{p}'})^{-1}(Y)\right)} \\
 & = & \tilde{g}_{s'}\left((\psi)_* \circ ((\varphi_s)_*|_{\tilde{p}})^{-1}(X),
(\psi)_* \circ ((\varphi_s)_*|_{\tilde{p}})^{-1}(Y)\right) \\
& = & \psi^* \tilde{g}_{s'} \left(((\varphi_s)_*|_{\tilde{p}})^{-1}(X),
((\varphi_s)_*|_{\tilde{p}})^{-1}(Y)\right) \\
& = & \tilde{g}_s \left(((\varphi_s)_*|_{\tilde{p}})^{-1}(X),
((\varphi_s)_*|_{\tilde{p}})^{-1}(Y)\right).
\ean

To show that $\varphi_s^* \bar{g} = \tilde{g}_s$, consider
\[
\varphi_s : \ B(0, \frac{9}{10}r) \ra U_s.
\]
 In particular, for any $\tilde{p} \in B(0, \frac{9}{10}r)$,
$B_{g_s}(\tilde{p},r_1) \subset B(0,r)$, and therefore $\tilde{p}$ can be used
to define the metric $\bar{g}$ at $\varphi_s (\tilde{p})$. It follows from the
definition that
\be  \label{mm}
\varphi_s^* \bar{g}  = \tilde{g}_s.
\ee
Finally, note that the smoothness of the metric $\bar{g}$ is an immediate
consequence of (\ref{mm}).
\qed

\sect{Embedding I}

We continue with the above manifold $(M, g)$. Let $\Omega = B(0, r) \subset
R^n$ with a pull back metric $\varphi_s^*g$. For convenience, this metric will
be denoted by $g$. It is easy  to see that $\|(\Omega, g)\|^h_{C^{0, \alpha},
r} \leq Q(r)$. We are going to construct an equivariant embedding of $(\Omega,
g)$ into $L^2(\Omega, g)$  by associating  to every point $p \in \Omega$ a
geometric function $f_p \in L^2(\Omega) \equiv L^2(\Omega, g)$, which depends
nicely on $p$. A natural choice seems to be the distance function measured from
$p$. Indeed, it is used by Abresch in \cite{a}. However, under our rather
weak assumptions on the metric it is impossible to have uniform control of the
second order derivative of the distance function, which is needed to ensure
that the pull-back metric induced by the embedding satisfies a sectional
curvature bound. In fact, one can not even expect differentiability of the
distance function in balls of uniform size.
Our substitute for the distance function is solutions of a
canonical geometric partial differential equation.
Those solution functions have the crucial equivariance (
like the distance functions) and enjoy better regularity.
Many choices of ``canonical" PDE solutions are possible, e.g. in \cite{a}
Green's function is suggested. But Green's function is inconvenient because of
its singularity.  We shall employ a very
simple  and nicely-behaved PDE.

Denote
\[
\Omega_1 = \Omega \setminus \cup_{q \in \partial \Omega} \overline{B_g(q,
i_0)},
\]
where $i_0 =\frac{r}{10}$. ($B^g(q, \cdot)$ denotes the closed geodesic ball of
center $q$ and radius $\cdot$ measured in $g$.)
Then for $s \in \Omega_1$, let $h_s  \in L^{1,2}_0 (\Omega)$ be the unique weak
solution of the following Dirichlet boundary value problem:
\be  \label{h}
\left\{ \begin{array}{rcll} \Delta h_s  & = & -1 & \mbox{in} \ B_g(s,i_0) \\
h_s   & \equiv & 0  & \mbox{on} \ \partial B_g(s,i_0).
\end{array}  \right.
\ee
Here the Laplace operator is defined with respect to the metric $g$. The
function $h_s$ will be extended to be zero outside the geodesic ball.

Since the harmonic $C^{0,\alpha}$-norm  of $(\Omega, g)$ is uniformly bounded,
it is easy to see that a uniform Poincare inequality holds on the balls
$B_g(s,i_0)$ with dependence on $i_0$. A simple integration argument then
yields a uniform estimate of the Sobolev norm of $h_s$. Uniform interior $C^{2,
\alpha}$ estimates then follow readily, because in harmonic coordinates the
Laplace operator takes the form $\Delta = g^{ij} \partial_i \partial_j$.    We
also have a uniform $L^{\infty}$ estimate up to boundary, but it seems
impossible to obtain better estimate up to boundary because the control of the
geometry of the boundary is very weak. At a first glance this appears to
threaten to destroy the embedding scheme. Fortunately we have a way to get
around it. On the other hand, we can not obtain control of the dependence of
$h_s$ on the center
$s$. To remedy this, we shall take a suitable average of $h_s$ over $s$. The
resulting new family of functions will depend nicely on the center.

Now let us state a few basic properties of the functions $h_s$ in the following
proposition, which will be proved at the end of this section.  Here, as before,
we work under the assumption $\|(\Omega, g)\|^h_{C^{0, \alpha}, r} \leq Q(r)$.

\begin{prop}  \label{bh}
 Let $\bar{h}_s(p)$ be the solution of equation (\ref{h}) with respect to the
canonical Euclidean metric $g_0$ on the Euclidean ball $B_{g_0}(s,i_0)$. Then
 for any $\epsilon > 0$ and fixed $0<R<1$, there is an $r_0 = r_0
(\epsilon,R,Q) > 0$ such that if $i_0 \leq r_0$,
\ba
|h_s(p) - \bar{h}_s(p)| & < & \epsilon i_0^2, \label{c1} \\
| \frac{\partial}{\partial p} h_s(p) - \frac{\partial}{\partial p}
\bar{h}_s(p)| & < &\epsilon i_0 \label{c2}
\ea
for all $s$ and all $p$ with $d_{g_0}(s,p) \leq Ri_0$. It will follow from the
proof that $B_{g_0}(s,Ri_0) \subset B_g(s, i_0)$ so that these estimate make
sense. Also
 \be  \label{ub}
| \frac{\partial^2}{\partial p^2} h_s(p) |\leq  C(n, Q,R), \ \ |\frac{1}{i_0}
\frac{\partial}{\partial p} h_s(p) |\leq  C(n, Q,R).
\ee
\end{prop}

Note that
\be
\bar{h}_s(p) = \frac{1}{2n} (i_0^2 - d_{g_0}^2(s,p)).
\ee
 Therefore $\frac{2n}{i_0^2} \bar{h}_s(p) \leq \frac{1}{5}$ when $d_{g_0}(s,p)
\geq \sqrt{\frac{4}{5}}i_0$. Choosing $R  = \frac{10}{11}$ in
Proposition~\ref{bh}, we have  $\frac{2n}{i_0^2} h_s(p) < \frac{1}{4}$ when $
\sqrt{\frac{4}{5}}i_0 \leq d_{g_0}(s,p) \leq  \frac{10}{11}i_0$ and $i_0$ is
sufficiently small.

Let $\beta = \beta_n \in C^\infty_0 ([0,\infty))$ be the cut off function.
\[
\beta_n(t) = \left\{ \begin{array}{ll} 0 & \mbox{if} \  0 \leq t \leq
\frac{1}{4} \\
  B_n & \mbox{if} \ t \geq \frac{1}{2}
\end{array} \right.,
\]
where
$B_n$ is a constant which will be determined later.

Then $\beta \left( \frac{2n}{i_0^2}h_s(p)\right) = 0$ near the sphere
$d_{g}(s,p) = \frac{9}{10}i_0$ for all $i_0$ small. (Note that $d_{g}$
converges to $d_{g_0}$ when $i_0 \ra 0$.) We define a new function which is
$\beta \left( \frac{2n}{i_0^2}h_s(p)\right)$ restricted to the ball $B(s,
\frac{9}{10}i_0)$
and identically zero outside. For simplicity we still denote this new function
by $\beta \left( \frac{2n}{i_0^2}h_s(p)\right)$.
 As mentioned before, we have no control of the dependence of $h_s$ on the
center $s$. The said average function is given as follows
\be
f_p(q) = \int_\Omega \beta \left( \frac{2n}{i_0^2}h_s(p)\right) \beta \left(
\frac{2n}{i_0^2}h_s (q) \right)\, ds.
\ee
Note that $f_p(q)$ is symmetric in $p$ and $q$ and is $C^{2,\alpha}$ uniformly
bounded in both variables.

Now we define the embedding
\ban
F: \ \Omega_1 & \ra & L^2(\Omega, g)\\
  p & \ra & i_0^{-\frac{3}{2}n+1} f_p(q)
\ean
Note that
\ba
d_{v_p}F : & q  \longmapsto & 2ni_0^{-\frac{3}{2}n}\int_\Omega \beta'\left(
\frac{2n}{i_0^2}h_s(p)\right)  \lp \frac{1}{i_0}\nabla_{v_p} h_s(p), v_p \rp
\beta \left( \frac{2n}{i_0^2}h_s (q) \right)\, ds, \label{dff}\\
\nabla^2_{v_p,w_p}F: & q \longmapsto & 4n^2i_0^{-\frac{3}{2}n-1}\int_\Omega
\beta''\left( \frac{2n}{i_0^2}h_s(p)\right)  \lp \frac{1}{i_0}\nabla h_s(p),
v_p \rp  \lp \frac{1}{i_0}\nabla h_s(p), w_p \rp \beta \left(
\frac{2n}{i_0^2}h_s (q) \right)\, ds  \nonumber \\
& & + 2ni_0^{-\frac{3}{2}n-1}\int_\Omega \beta'\left( \frac{n}{i_0^2}\right)
\nabla^2_{v_p,w_p} h_s(p) (v_p,w_p) \beta \left( \frac{n}{i_0^2}h_s (q)
\right)\, ds.  \label{dff2}
\ea

We first show that when $\Omega$ is an Euclidean domain, we can normalize
$\beta$ so that $F$ is an isometric imbedding. In this case the imbedding
function
\[
\bar{f}_p(q) = \int_\Omega \beta \left(1-\frac{d_{g_0}^2(p,s)}{i^2_0}\right)
\beta \left(1-\frac{d_{g_0}^2(q,s)}{i^2_0}\right)\, ds.\]
By the symmetry of the integration domain, $B(p, \frac{\sqrt{3}i_0}{2}) \cap
B(q, \frac{\sqrt{3}i_0}{2})$, and the integrand, $\bar{f}_p(q)$ depends only on
$d_{g_0}(p,q)$ (and $\beta$). We can write $\bar{f}_p(q) = i^n_0 \tilde{f} (n,
\frac{1}{i_0}d_{g_0}(p,q))$,
$ d_{v_p}\bar{F}(q) =  i_0^{-n/2}\tilde{f}' (n, \frac{1}{i_0}d_{g_0}(p,q)) \lp
\nabla d_{g_0}(p,q), v_p \rp$. Then
\ban
\| d_{v_p}\bar{F}\|^2_{L^2(\Omega)} & = &  i_0^{-n}\int_{B(p, 2i_0)}
\tilde{f}'^2 (n, \frac{1}{i_0}d_{g_0}(p,q)) \lp \nabla d_{g_0}(p,q), v_p \rp^2
dq \\
& = &  i_0^{-n}\int_0^{2i_0} r^{n-1} \int_{S^{n-1}}  \tilde{f}'^2 (n,
\frac{1}{i_0}r) \lp \xi, v \rp^2 d\xi dr \\
& = &  i_0^{-n}\frac{\vol (S^{n-1})|v|^2}{n}\int_0^{2i_0} r^{n-1}\tilde{f}'^2
(n,  \frac{1}{i_0}r) dr \\
& =& \frac{\vol (S^{n-1})|v|^2}{n}\int_0^{2}r^{n-1}\tilde{f}'^2 (n,r) dr.
\ean
Choose $B_n$ in the defintion of $\beta$ so that $\frac{\vol
(S^{n-1})}{n}\int_0^{2}r^{n-1}\tilde{f}'^2 (n,r) dr = 1$. Then we will have
achieved the following.
\begin{lem}
$\bar{F}$ is an isometric embedding.
\end{lem}

With the above choice of $\beta$ we will show that $F$ is an almost isometric
embedding when $\Omega$ is not necessarily an Euclidean domain and the second
derivative of $F$ is also uniformly bounded. More precisely we have
\begin{prop} \label{e-df}
For any given $\epsilon_0 >0$, there exists an $r_0>0$ such that
\be  \label{df1}
(1+ \epsilon_0)^{-2} |v|^2 \leq \|d_{v_p}F\|^2_{L^2(\Omega)} \leq (1+
\epsilon_0)^{2} |v|^2,
\ee
for all $v\in T_p\Omega_1$ and $0 < i_0 \leq r_0$. And
\be  \label{df2}
\|\nabla^2_{v_p,w_p}F\|^2 _{L^2(\Omega)} \leq C(n,\alpha,Q) i_0^{-2}|v|^2 \cdot
|w|^2.
\ee
\end{prop}
\Pf By definition
\[
\|d_{v_p}F\|^2_{L^2(\Omega)} = \int_{B(p,2i_0)} |d_{v_p}F(q)|^2 dq.
\]
Now the volume element of the metric $g$ is comparable with Euclidean one.
Namely
\be  \label{vv}
e^{-Q(i_0)} \vol_{R^n} \leq \vol_g \leq e^{Q(i_0)} \vol_{R^n}.
\ee
Therefore it suffices to prove that $|d_{v_p}F(q)|$ is close to
$|d_{v_p}\bar{F}(q)|$ when $i_0$ is small, which follows from  (\ref{dff}) and
Proposition~\ref{bh}.

(\ref{df2}) also follows from (\ref{dff2}), (\ref{vv}) and
Proposition~\ref{bh}.
\qed  \\

{\em Proof of} Proposition~\ref{bh}. First we introduce some new functions.
Let $\tilde{h}_{s, i_0}(p)$ be the solutions of (\ref{h}) on
$B_{i_0^{-2}g}(s,1)$ with  respect to the scaled metrics $i_0^{-2}g$, and
$\bar{\tilde{h}}_s(p)$ the solution of (\ref{h}) with respect to the Euclidean
metric on the Euclidean ball $B(s,1)$. Then
\be  \label{hh}
\tilde{h}_{s,i_0}(p)= i_0^{-2} h_s(p), \ \ \bar{\tilde{h}}_s(p) = i_0^{-2}
\bar{h}_s(p).
\ee
 (Here the variables $s,p$ are in domain with different metrics for different
functions.)  Since
 \be \label{normb}
\|(B_{i_0^{-2}g}(s,1), i_0^{-2}g)\|_{C^{0,\alpha},1} = \|(B_g(s,i_0),
g)\|_{C^{0,\alpha},i_0} \leq Q(i_0) \leq Q(1), \ \mbox{for}\ i_0 \leq 1,
\ee
for the same reasons as mentioned before for $h_s$, we have the following
estimates
\be  \label{el}
\| \tilde{h}_{s,i_0}\|_{L_0^{1,2}(B_g(s,1))} \leq C(n,Q(1))
\ee
and
for any fixed $0<R<1$,
 \be  \label{ell}
\| \tilde{h}_{s,i_0} \|_{C^{2,\alpha}(B_g(s,R))}  \leq C(n, Q(1),R).
\ee

On the other hand, the hypothesis $\|(B_{i_0^{-2}g}(s,1),
i_0^{-2}g)\|_{C^{0,\alpha},1} \leq Q(i_0)$ implies that
\[
e^{-Q(i_0)} d_{i_0^{-2}g_0}(p,s) \leq d_{i_0^{-2}g}(p,s) \leq e^{Q(i_0)}
d_{i_0^{-2}g_0}(p,s).
\]
Therefore
\[
B_{i_0^{-2}g_0} (s, e^{-Q(i_0)}R) \subset B_{i_0^{-2}g} (s, R) \subset
B_{i_0^{-2}g_0}(s, e^{Q(i_0)}R).
\]

{}From these estimates and the uniqueness of the weak solution $\bar{\tilde
h}_s$ it is easy to deduce the following: for each sequence of centers $s_k$
converging to some center $s_0$ and each sequence $i_0(k)$ converging to zero,
the corresponding rescaled solutions $\tilde h_{s_k, i_0(k)}$ converge weakly
to $\bar{\tilde h}_{s_0}$. Moreover, by the Arzela-Ascoli theorem,  they also
converge uniformly in $C^1$ on proper compact subsets of $B_{g_0}(s_0, 1)$.
This convergence fact along with the smooth dependence of $\bar{\tilde h }_s$
on $s$ then imply that the $\tilde h_{s, i_0}$ converge uniformly with respect
to $s$ in $C^1$ on proper compact subsets of $B_{g_0}(s, 1)$ as $i_0$ goes to
zero. Consequently, for a fixed $R \in (0,1)$, given any $\epsilon >0$,
there is an $r_0 >0$ such that  for all $s$ and $p$ with $d_{g_0}(p,s) < R$, if
$i_0 \leq r_0$, then
\be
 |\tilde{h}_s(p) -\bar{\tilde{h}}_s(p)|  < \epsilon. \label{pc1}
\ee
Similarly,
\be
| \frac{\partial}{\partial p} \tilde{h}_s(p) -\frac{\partial}{\partial p}
\bar{\tilde{h}}_s(p)|  < \epsilon. \label{pc2}
\ee

Hence for all $s$ and all $p$ with $d_{g_0}(s,p) \leq Ri_0$,
\ban
|h_s(p) - \bar{h}_s(p)| & < & \epsilon i_0^2,  \\
| \frac{\partial}{\partial p} h_s(p) - \frac{\partial}{\partial p}
\bar{h}_s(p)| & < &\epsilon i_0.
\ean

(\ref{ub}) just follows from (\ref{ell}) and (\ref{hh}).
\qed

\sect{Embedding II}
In this section we study the geometry of $F(\Omega_1)$ as a submanifold in
$L^2(\Omega)$. We will prove, among other things, two important properties of
$F(\Omega_1)$. That is, the induced metric of $F(\Omega_1)$ has uniformly
bounded sectional curvature and the embedding $F$ is equivariant.

The geometry of $F(\Omega_1)$ is completely determined by the second
fundamental form of its embedding into $L^2(\Omega)$, which in turn can be
described by the family of orthogonal projections.
$P(y):\ L^2(\Omega) \ra T_yF(\Omega_1) \subset L^2(\Omega), y \in F(\Omega^1)$.
We have
\begin{lem}  \label{4.1}
 The sectional curvature of $F(\Omega_1)$ is given by the following formula:
\be \label{curv}
R(z_1, z_2) z_3 = [d_{z_1} P, d_{z_2}P] z_3, \ \ z_1,z_2,z_3 \in T_z
F(\Omega_1).
\ee
\end{lem}
\Pf Since $P^2 =P$,  one has
\be  \label{p2}
(d_{z_1}P)P + P (d_{z_1}P) = d_{z_1}P.
\ee
Let $\nabla$ be the connection on $F(\Omega_1)$ and $d_{z_1}$ the directional
derivative on the $L^2$ space. Then
\ban
\nabla_{z_1} z_2 & = & P(d_{z_1}z_2) \\
& = & d_{z_1}z_2 - (1- P)(d_{z_1}(Pz_2)) \\
& = & d_{z_1}z_2 - (1-P) \left[ (d_{z_1}P)z_2 + P (d_{z_1}z_2) \right] \\
& = & d_{z_1}z_2 - (1-P) (d_{z_1}P) (Pz_2) \\
& = & d_{z_1}z_2 - (d_{z_1}P)z_2.
\ean
Here we have used (\ref{p2}) in the last equation.
Therefore
\be   \label{p1}
\nabla_{z_1}  = d_{z_1} - (d_{z_1}P).
\ee
Now formula \ref{curv} follows from (\ref{p1}) and the definition of the
curvature tensor.
\qed

\begin{prop}  \label{4.2}
Let $\alpha_0 = (1+\epsilon_0)^2 C(n,\alpha,Q)$. Here $\epsilon_0, r_0, C$ are
the same constants as in Proposition~\ref{e-df}. Then for all $0 < i_0 < r_0$,
\be  \label{ep}
\| d_{\dot{y}}P\|_{op} \leq \alpha_0 i_0^{-1} \| \dot{y} \|, \ee
\end{prop}
\Pf Since $\left( 1-P(F(p)) \right) d_{w_p} F = 0$,
\[
d_{d_{v_p}F} P \cdot d_{w_p} F = \left( 1-P(F(p)) \right) \nabla^2_{v_p,w_p}F.
\]
By (\ref{df1}) and (\ref{df2}), $\| d_{\dot{y}}P\|_{op} \leq \alpha_0 i_0^{-1}
\| \dot{y} \|$.
\qed \\

Therefore the metric $\tilde{g} = F^* g_{L^2}$,  the metric on $\Omega_1$
obtained by pulling back the $L^2$ metric, has bounded sectional curvatures.

To prove the equivariance, we first note:
\begin{lem}  \label{uh}
Let $h_s(p)$ be the function defined in (\ref{h}), and let  $\psi: \Omega \ra
\Omega'$ be an isometric embedding. Then
\be
h_{\psi (s)} (\psi (p)) = h_s(p).
\ee
\end{lem}
\Pf Since equation (\ref{h}) is invariant under isometry, this follows from the
uniqueness of solutions to (\ref{h}).
\qed \\

Let $(\Omega,g)$ be as before and $F: \Omega_1 \ra L^2(\Omega)$ the embedding
defined in \S 3. With the above lemma, we can now prove
\begin{prop}  \label{4.4}
If $\psi : (\Omega,g) \ra (\Omega',g')$ is an isometric embedding, then $$\psi
: (\Omega, \tilde{g}) \ra (\Omega', \tilde{g}')$$ is also an isometric
embedding.
\end{prop}
\Pf First, we assume $\psi$ is actually an isometry. Then
\[
F \circ \psi (p) = f_{\psi (p)},
\]
where the function
\ban
f_{\psi (p)} (q) & = &  \int_\Omega \beta \left( \frac{2n}{i_0^2}h_s(\psi
(p))\right) \beta \left( \frac{2n}{i_0^2}h_s (q) \right)\, ds.
\\
& = &  \int_\Omega \beta \left( \frac{2n}{i_0^2}h_{\psi^{-1}(s)}(p)\right)
\beta \left( \frac{2n}{i_0^2}h_{\psi^{-1}(s)} (\psi^{-1}(q)) \right)\, ds.
\ean
Here we have used Lemma~\ref{uh}. Since $\psi$ is an isometry, a change of
coordinates yields
\[
f_{\psi (p)}(q) = f_p (\psi^{-1} (q)).
\]
It follows then that
\[
F \circ \psi = (\psi^{-1})^* \circ F,
\]
where we have denoted by $(\psi^{-1})^*$ the map on $L^2(\Omega)$ induced by
$\psi^{-1}$. Therefore
\[
\psi^* F^*g_{L^2} = (F \circ \psi)^* g_{L^2} = F^* ((\psi^{-1})^*)^* g_{L^2} =
F^*g_{L^2}.
\]
This proves the equivariance when $\psi$ is an isometry. Since $\Omega,
\Omega'$ are both domains of $R^n$, the general statement follows by applying
the above to $\psi:\ \Omega \ra \psi (\Omega)$.
\qed \\

{\it Proof of}  Theorem 1.1.  This theorem is a consequence of Proposition 2.2,
Lemma 4.1, Proposition 4.2 and Proposition 3.3.
\qed

\sect{Proof of Theorem 1.2}
We consider $(M,g)  \in {\cal M}(n,p,Q)$ and $\Omega = B(0, r) \subset R^n$
with the pull-back metric $\varphi_s^*g$, where $\varphi_s$ is a coordinate
map. For convenience, we shall again  denote the pull-back    metric by $g$. We
have the inequality $\|(\Omega, g)\|_{L^{1,p}, r} \leq Q(r)$.

We employ the same embedding of $(\Omega, g)$ as before. Thus we use the same
functions $h_s \in L_0^{1,2}(\Omega)$, as given by (3.1). But the estimates for
$h_s$ are different now. By the $L^p$ elliptic theory we have uniform $L^{2,p}$
estimates for $h_s$ in the interior. A result similar to Proposition 3.1 then
holds, namely we have the estimates (3.2), (3.3)  and the second one in (3.4),
while the first one in (3.4) is replaced by an estimate on the $L^p$-norm of
the second order derivative.  Now it  is clear that our smoothing process
produces a metric $\bar g$ on $M$ such that
the lifted metrics $\varphi_s^* \bar g$ on $B(0, r/2)$ have uniformly bounded
Sobolev constant and  uniformly $L^p$-bounded sectional curvatures.

Next we apply the Ricci flow to deform $\bar g$. For this purpose, we assume
that $M$ is closed. We appeal to the arguments in \cite{dwy}. There, manifolds
with a pointwise bound on Ricci curvature and a conjugate radius bound are
treated.
These conditions are used to show that controlled harmonic coordinates exist
on lifted local  patches, where the lifting is given by the exponential map.
In these coordinates, the Ricci curvature bound then implies an $L^p$-bound on
sectional curvatures. In our situation, we do have an $L^p$-bound
on sectional curvatures. But the Ricci curvature bound is also used in several
other places  in \cite{dwy}. Since this bound is not available here, we need to
modify the arguments in \cite{dwy}.

In the key Proposition 3.1 (uniform short time
existence of the Ricci flow with a priori control) in \cite{dwy}, we drop the
estimate (3.4) on Ricci curvature. (Note that the conclusion of the proposition
without (3.4) still suffices for our purpose.)
We  claim that the proposition then holds in our new situation. For
convenience, we shall call this proposition the "Key Proposition". In
\cite{dwy}, the proof of Key Proposition is based on four lemmata: Lemmata 3.2,
3.3, 3.4 and 3.5.
Now let's take  a look at these lemmata in our new situation. Lemma 3.2
(pointwise estimate for Riemann curvature tensor) holds without change.  The
proof of Lemma 3.3 ($L^p$-estimate for Riemann curvature tensor) depends on a
{\it covering  estimate}, i.e. an estimate  for certain covering number and
multiplicity regarding geodesic balls in the lifted patches. In \cite{dwy}, the
estimate comes from the Bishop-Gromov covering argument, which depends on a
pointwise lower bound for Ricci curvature.
Now we do not have  such a bound.   But we still have a covering estimate,
which follows from the properties of our coordinates and the basic control
over the pull-back metrics as given by (the modified versions of) Propositions
3.1  and
3.2 (in the present paper).

 Another ingredient in the proof of Lemma 3.3 is
an isometry  correspondance between geodesic balls on different lifted patches.
Now it is given by Lemma~\ref{hom}. (The proof of this correspondence given in
\cite{dwy} does not work here.) Thus Lemma 3.3 also holds.  Since we have
dropped the estimate (3.4) about Ricci curvature in Key Proposition, Lemma 3.4
is no longer needed. Finally, note that the proof for Lemma 3.5 (estimate of
the Sobolev constant) in \cite{dwy} goes by computing the change rate of the
Sobolev constant along the flow. In \cite{dwy}, this rate is controlled by
a uniform bound on Ricci curvature,  which is not valid here.
However, Lemma 3.2 contains an estimate for Ricci curvature at positive time
$t$, namely it is
dominated by a constant times $t^{-1/2}$. Since the function $t^{-1/2}$ is
integrable at $0$, it is clear  that  the change rate of the Sobolev constant
is still under control without a uniform Ricci curvature bound, and hence Lemma
3.5 carries over.  (An alternative way of handling the Sobolev constant is to
apply Yang's estimate for it in \cite{y2}, which uses only an $L^p$-bound on
Ricci curvature and a positive lower bound on (local) volume. But that is more
involved. )

We leave to the reader to formulate precisely the  independent result implied
by the above proof about short time existence of the Ricci flow.

P. Petersen

Department of Mathematics,
UCLA
Los Angles, CA  90024

 petersen@math.ucla.edu

G. Wei \& R. Ye

Department of Mathematics,
 University of California,
Santa Barbara, CA 93106

 wei@math.ucsb.edu \ \ \
yer@math.ucsb.edu

R.Ye

Mathematics Institue,
Bochum University,
44780 Bochum

\end{document}